\documentclass[11pt]{amsart}
\usepackage[a4paper,margin=27mm]{geometry}
\usepackage[T1]{fontenc}
\usepackage{lmodern}
\usepackage{amsmath,amssymb,mathtools}
\usepackage{microtype}
\usepackage{graphicx}
\usepackage{tikz}
\usetikzlibrary{arrows.meta,decorations.pathreplacing}
\usepackage{placeins}
\usepackage{needspace}
\usepackage[colorlinks=true,linkcolor=blue,citecolor=blue,urlcolor=blue]{hyperref}

\newtheorem{theorem}{Theorem}[section]
\newtheorem{lemma}[theorem]{Lemma}

\newtheorem{corollary}[theorem]{Corollary}
\theoremstyle{definition}
\newtheorem{definition}[theorem]{Definition}
\newtheorem*{definition*}{Definition}
\newtheorem{assumption}[theorem]{Assumption}
\theoremstyle{remark}
\newtheorem{remark}[theorem]{Remark}
\numberwithin{equation}{section}
\allowdisplaybreaks[1]
\title[Generic linear momentum escape]{Generic Linear Momentum Escape for Twist Maps}
\author{Jianxing Du}
\address{School of Mathematical Sciences, Beijing Normal University,
No. 19, XinJieKouWai St., HaiDian District, Beijing 100875, P. R. China}
\email{jianxingdu@mail.bnu.edu.cn, jxdu000@gmail.com}
\hypersetup{pdfauthor={Jianxing Du},pdftitle={Generic Linear Momentum Escape for Twist Maps}}
\date{}
\keywords{Exact symplectic twist maps, generic escape, stability,
variational methods, invariant curves}

\begin{document}
\begin{abstract}
We prove that linear momentum escape is generic in the $C^{1,\alpha}$
topology for smooth potential perturbations of exact symplectic twist
maps, for every $0<\alpha<1$.
We assume a $C^2$ generating function, a uniform twist condition, and
uniform upper bounds on $D_{11}S$ and $D_{22}S$.
For an open dense set of potentials, every sufficiently small momentum
drift is realized by an orbit.
This genericity result does not extend in general to the $C^2$ topology.
We give examples with no nonzero linear drift under $C^2$-small
perturbations in every dimension.
In dimension one, a second example has bounded momentum along every
orbit under $C^2$-small perturbations.
\end{abstract}
\maketitle

\section{Introduction}\label{sec:introduction}

We study momentum escape for exact symplectic twist maps on
$\mathbb T^d\times\mathbb R^d$. We fix a generating function $S$ and
perturb it by a smooth periodic potential. We prove that linear
momentum escape is generic in the relative $C^{1,\alpha}$ topology
on smooth potentials, for every $0<\alpha<1$.
Here \emph{generic} means that the property holds on an open dense set.
For every potential in this set, the perturbed map has an orbit with
each sufficiently small linear momentum drift.

We first construct one small potential that realizes a whole ball of
drifts. We then carry out the construction near any smooth potential
and prove that the drift ball persists under small $C^0$ perturbations
of the potential. We also give examples showing that this genericity
result does not extend in general to the $C^2$ topology.

The question is motivated by Mather's conjecture on generic escape for
time-periodic Lagrangian systems~\cite[Section~13]{Mather1993}.
Our results concern a discrete system and weaker topologies on the
potentials. We explain the relation in Section~\ref{sec:mather-conjecture}.

\subsection{Setting and linear momentum drift}

We write $\mathbb T^d=\mathbb R^d/\mathbb Z^d$ and denote by
$\pi_p:\mathbb T^d\times\mathbb R^d\to\mathbb R^d$ the momentum
projection. Vector norms are Euclidean, and matrix norms are the
associated operator norms. We write $D_iS$ for the gradient in the
$i$th endpoint and $D_{ij}S$ for the corresponding Hessian block.
Thus $D_{21}S=D_{12}S^{\mathsf T}$.
For a smooth function $U$ on the flat torus,
we use
\[
 \|U\|_{C^r}=\max_{0\leq j\leq r}\sup_x\|D^jU(x)\|,
 \qquad
 \|U\|_{C^{1,\alpha}}=\|U\|_{C^1}+[DU]_{\alpha},
\]
where the H\"older seminorm is computed using the flat torus distance.
Periodic functions are identified with their lifts to $\mathbb R^d$.

We consider generating functions $S$ satisfying Assumption~\ref{ass:generating-function}
in Section~\ref{subsec:twist-maps}.

For $V\in C^\infty(\mathbb T^d)$, let $F_V$ be the map generated by
\[
 S_V(x,y)=S(x,y)-V(x).
\]
Its lift is defined by
\begin{equation}\label{eq:generating-equations}
 \widetilde F_V(x,p)=(y,P)
 \quad\Longleftrightarrow\quad
 \begin{cases}
 p=-D_1S(x,y)+DV(x),\\
 P=D_2S(x,y).
 \end{cases}
\end{equation}
These equations define a global $C^1$ diffeomorphism whose quotient is
an exact symplectic twist map. We recall these facts in
Section~\ref{subsec:twist-maps}.

\begin{definition*}[Momentum drift and escape]
Let $z\in\mathbb T^d\times\mathbb R^d$ and write
$p_n=\pi_pF_V^n(z)$ for $n\in\mathbb Z$.
In this paper, the orbit has \emph{linear momentum drift}
$p^*\in\mathbb R^d$ if
\[
 \sup_{n\in\mathbb Z}|p_n-np^*|<\infty.
\]
When $p^*\neq0$, we call this \emph{linear momentum escape}.
The orbit is \emph{escaping} if $|p_n|\to\infty$ as
$n\to+\infty$.
\end{definition*}

A nonzero drift gives escape in both time directions. Zero drift gives
bounded momentum. Escape alone does not prescribe a growth
rate or a limiting direction.

\begin{theorem}\label{thm:acceleration}
Let $d\geq1$, and let $S$ satisfy Assumption~\ref{ass:generating-function}.
For every $\varepsilon>0$ and every $0<\alpha<1$, there exist
$V\in C^\infty(\mathbb T^d)$ and $r>0$ such that
\[
 \|V\|_{C^{1,\alpha}}<\varepsilon
\]
and, for each $p^*\in\mathbb R^d$ with $|p^*|\leq r$, there exists
$z_{p^*}\in\mathbb T^d\times\mathbb R^d$ satisfying
\begin{equation}\label{eq:bounded-drift-error}
 \sup_{n\in\mathbb Z}
 \big|\pi_pF_V^n(z_{p^*})-np^*\big|<\infty.
\end{equation}
Thus,
\[
 \lim_{n\to\infty}\frac{\pi_pF_V^n(z_{p^*})}{n}=p^*.
\]
\end{theorem}

The potential $V$ is chosen before $p^*$. Thus one potential realizes
the whole drift ball. The proof bounds the supremum
in~\eqref{eq:bounded-drift-error} by $2\sqrt d\,K^2/\lambda$,
uniformly for $|p^*|\leq r$.

Accelerating orbits of twist maps have been studied through shear and
vertical rotation numbers. Slijep\v{c}evi\'{c}~\cite[Theorem~3.7]{Slijepcevic1999}
constructed orbits with independently prescribed forward and backward
shear-rotation numbers in the shear-rotation interval of a twist
diffeomorphism of $\mathbb T^2$.
Taking both numbers equal to $\rho$ gives $p_n=n\rho+o(|n|)$ as
$|n|\to\infty$. The theorem does not give a common bound for
$|p_n-n\rho|$ over the whole interval.

He also stated a large-amplitude result for
$S(x,y)+kP(x)$~\cite[Proposition~4.3]{Slijepcevic1999}.
Here $P$ is a fixed $C^2$ periodic potential whose critical points
are nondegenerate. For every positive integer $m$, the shear-rotation
interval contains $[-m,m]$ for all sufficiently large $k$.
Theorem~\ref{thm:acceleration} instead gives arbitrarily small
$C^{1,\alpha}$ potentials.

Addas-Zanata~\cite{AddasZanata2005} studied vertical rotation intervals
and their behavior under perturbations. The vertical rotation number
of an orbit $(x_n,p_n)$ is $\lim_{n\to\infty}(p_n-p_0)/n$, when this
limit exists. Periodic accelerator orbits already satisfy
$p_n=n\rho+O(1)$. This bound also follows from the compact sets at
rational endpoints recalled in~\cite[Theorem~3.6]{AddasZanata2025}.
These works assume an additional periodicity in momentum, which
allows the maps to descend to $\mathbb T^2$.

Small potential perturbations can also destroy invariant circles.
For the quadratic generating function $S(x,y)=\tfrac12(y-x)^2$,
Wang~\cite[Section~2]{Wang2014} gives analytic potential perturbations
that are arbitrarily small in $C^{2-\delta}$ and destroy all essential
invariant circles. Taking $\delta=1-\alpha$ gives the $C^{1,\alpha}$
norm used here. Thus the destruction of invariant circles at this
regularity is already known. Here the norm is that of the potential,
rather than that of the map.

\subsection{Generic linear momentum escape}

For each fixed $0<\alpha<1$, we give $C^\infty(\mathbb T^d)$
the relative $C^{1,\alpha}$ topology.
The next theorem strengthens Theorem~\ref{thm:acceleration} in two ways.
The constructed potential can approximate any given smooth potential
in the $C^{1,\alpha}$ norm. The whole drift ball persists under small
$C^0$ perturbations of the potential.

\begin{theorem}
\label{thm:robust-drift}
Let $S$ satisfy Assumption~\ref{ass:generating-function}.
For every $U\in C^\infty(\mathbb T^d)$, every
$\varepsilon>0$, and every $0<\alpha<1$, there exist
$V\in C^\infty(\mathbb T^d)$ and constants $\delta,r,C>0$ such that
\begin{equation}\label{eq:robust-approximation}
 \|V-U\|_{C^{1,\alpha}}<\varepsilon.
\end{equation}
For every $E\in C^\infty(\mathbb T^d)$ with
\begin{equation}\label{eq:nearby-potential}
 \|E\|_\infty<\delta
\end{equation}
and every $p^*\in\mathbb R^d$ with $|p^*|\leq r$, there is a point
$z_{E,p^*}\in\mathbb T^d\times\mathbb R^d$ such that
\begin{equation}\label{eq:persistent-drift}
 \sup_{n\in\mathbb Z}
 \big|\pi_pF_{V+E}^{\,n}(z_{E,p^*})-np^*\big|
 \leq C.
\end{equation}
The radius $r$ and the bound $C$ are independent of $E$.
\end{theorem}

The initial point may depend on $E$ and $p^*$.
The constants $\delta$, $r$, and $C$ are chosen before $E$ and $p^*$.

More recently, Addas-Zanata~\cite[Theorems~2.5 and~2.6]{AddasZanata2025}
proved that vertical rotation intervals are locally constant on
$C^r$-open dense sets of torus twist maps, for every integer $r\geq1$.
Their endpoints are rational. In the area-preserving class, these
intervals have nonempty interior. The local constancy holds under
$C^0$ perturbations of the map.

The density proof in~\cite[Section~4.2]{AddasZanata2025}
uses nonzero vertical translations. These do not preserve exactness.
Here we fix $S$ and vary only the periodic potential $V$, so every
$F_V$ remains exact symplectic. Thus density in the full
area-preserving class does not imply density in this potential family.
Our results apply to every $S$ satisfying
Assumption~\ref{ass:generating-function}, in every dimension and
without periodicity in momentum.
Theorem~\ref{thm:robust-drift} preserves the same drift ball and a
common bound under $C^0$-small changes of the potential.
No bound on derivatives of $E$ is required. Since the map depends on
the gradient of the potential, this persistence does not follow from
stability under $C^0$ perturbations of the map.

\begin{corollary}
\label{cor:open-dense}
Let $S$ satisfy Assumption~\ref{ass:generating-function} and let $0<\alpha<1$.
There is an open dense subset of $C^\infty(\mathbb T^d)$ in the relative
$C^{1,\alpha}$ topology. For each $V$ in this subset, there is $r(V)>0$
such that, for every $p^*\in\mathbb R^d$ with $|p^*|\leq r(V)$, some
$z_{p^*}\in\mathbb T^d\times\mathbb R^d$ satisfies~\eqref{eq:bounded-drift-error}.
\end{corollary}

We prove the corollary in Section~\ref{sec:persistence}.

\subsection{Counterexamples for \texorpdfstring{$C^2$}{C2}-small perturbations}

The $C^{1,\alpha}$-generic conclusion cannot in general be strengthened
to a $C^2$-generic conclusion. The next example shows that nonzero
linear momentum drift can be absent on a $C^2$ neighborhood of the zero
potential in every dimension. A second example excludes all escape in
dimension one. Both examples use smooth generating functions with
bounded Hessians.

\begin{theorem}
\label{thm:slowdown}
For every $d\geq1$, there exist
$S\in C^\infty(\mathbb R^{2d})$ and $\varepsilon_0>0$ such that
$S$ satisfies Assumption~\ref{ass:generating-function},
$\|D^2S\|_\infty<\infty$, and, for every
$V\in C^\infty(\mathbb T^d)$ with $\|V\|_{C^2}\leq\varepsilon_0$
and every $z\in\mathbb T^d\times\mathbb R^d$,
\begin{equation}\label{eq:zero-liminf}
 \liminf_{n\to+\infty}\frac{|\pi_pF_V^n(z)|}{n}=0.
\end{equation}
\end{theorem}

For this $S$, let $V$ satisfy $\|V\|_{C^2}\leq\varepsilon_0$.
If $\pi_pF_V^n(z)/n$ converges as $n\to\infty$, its limit is zero.
Potentials with nonzero linear momentum drift are not
$C^2$-dense. Hence the $C^{1,\alpha}$ norm in
Theorem~\ref{thm:acceleration} cannot in general be replaced by the
$C^2$ norm, and Corollary~\ref{cor:open-dense} has no general $C^2$
version.
However, the estimate~\eqref{eq:zero-liminf} still allows sublinear escape.

In dimension one, a different example prevents escape altogether
throughout a $C^2$ neighborhood of the zero potential.

\begin{samepage}
\begin{theorem}\label{thm:bounded}
There exist $S\in C^\infty(\mathbb R^2)$ and $\varepsilon_0>0$
such that $S$ satisfies Assumption~\ref{ass:generating-function},
$\|D^2S\|_\infty<\infty$, and, for every
$V\in C^\infty(\mathbb T)$ with $\|V\|_{C^2}\leq\varepsilon_0$
and every $z\in\mathbb T\times\mathbb R$,
\begin{equation}\label{eq:bounded-momentum}
 \sup_{n\in\mathbb Z}|\pi_pF_V^n(z)|<\infty.
\end{equation}
\end{theorem}
\end{samepage}

The bound in~\eqref{eq:bounded-momentum} may depend on $V$ and $z$.
The two examples use different generating functions.

For $S(x,y)=\tfrac12(y-x)^2$, KAM theory gives an invariant essential
circle when the potential is sufficiently small in a higher
differentiability norm~\cite{KatokKrikorian2022,Russmann1970}.
Its integer translates in momentum bound every orbit. This already
excludes $C^\infty$-generic escape in dimension one.
Theorem~\ref{thm:bounded} constructs a different generating function
for which every orbit has bounded momentum while controlling only
the $C^2$ norm of the potential. Thus escape need not be
$C^2$-generic either.

\subsection{Relation to Mather's escape conjecture}
\label{sec:mather-conjecture}

In \emph{Variational construction of connecting orbits},
Mather~\cite{Mather1993} gives variational conditions for connecting
action-minimizing sets. He also constructs orbits that visit prescribed
sequences of such sets. The paper ends with an escape conjecture
\cite[Section~13, pp.~1383--1384]{Mather1993}.

Let $M$ be a smooth compact manifold of dimension at least two, and let
$L_0:TM\times\mathbb T\to\mathbb R$ be a smooth time-periodic Lagrangian.
Its velocity Hessian is assumed to be positive definite, its growth
superlinear, and its Euler--Lagrange flow complete.
Assume also that the flow remains complete for every
\[
 L_\psi(q,\dot q,t)=L_0(q,\dot q,t)+\psi(q,t),
 \qquad \psi\in C^\infty(M\times\mathbb T).
\]
Mather conjectures that a residual set of these potentials in the
$C^\infty$ topology gives a flow with an orbit that eventually leaves
every compact subset of $TM\times\mathbb T$.
The conjecture does not prescribe an escape rate.

Our construction uses finite action minimization, as in discrete
Aubry--Mather theory; see Garibaldi and
Thieullen~\cite{GaribaldiThieullen2011}.
Corollary~\ref{cor:open-dense} gives a discrete
$C^{1,\alpha}$-generic escape result with prescribed linear momentum
drift. It does not establish Mather's conjecture, which concerns the
$C^\infty$ topology and continuous-time flows. Our potential perturbation is a momentum kick
depending on position; see~\eqref{eq:potential-kick}.
We do not identify it with an arbitrary smooth time-periodic potential.

The counterexamples also do not disprove Mather's conjecture.
Theorem~\ref{thm:slowdown} rules out nonzero linear momentum drift but
allows sublinear escape. Theorem~\ref{thm:bounded} concerns dimension
one, which is excluded from the conjecture.

\subsection{Outline of the proofs}

For Theorem~\ref{thm:acceleration}, finite action minimization gives
orbits with prescribed momentum drift for an auxiliary generating
function with an affine force. We choose a small periodic potential
so that a periodic correction cancels this force along the minimizing
configurations. A finite difference estimate gives the stability in
Theorem~\ref{thm:robust-drift}.

For Theorem~\ref{thm:slowdown}, we make momentum increments small on
wide annuli and estimate the time needed to cross them.
For Theorem~\ref{thm:bounded}, invariant essential circles bound
every orbit.

Section~\ref{sec:preliminaries} gives the required lemmas.
Sections~\ref{sec:acceleration}--\ref{sec:bounded} prove the four theorems.

\section{Preliminaries}\label{sec:preliminaries}

\subsection{Twist maps and potential perturbations}\label{subsec:twist-maps}

The phase space $\mathbb T^d\times\mathbb R^d$ carries the symplectic form
\[
 \omega=\sum_{j=1}^d dx_j\wedge dp_j.
\]
A $C^1$ diffeomorphism $F$ is \emph{symplectic} if $F^*\omega=\omega$.
We introduce twist maps through their generating functions, following
Arnaud's description of globally positive diffeomorphisms~\cite{Arnaud2016}.

\begin{definition}\label{def:twist-map}
An \emph{exact symplectic twist map} is a symplectic
$C^1$ diffeomorphism
$F:\mathbb T^d\times\mathbb R^d\to\mathbb T^d\times\mathbb R^d$
that is homotopic to the identity and has a lift
$\widetilde F:\mathbb R^d\times\mathbb R^d\to\mathbb R^d\times\mathbb R^d$
with a function $S\in C^2(\mathbb R^d\times\mathbb R^d)$ satisfying:
\begin{enumerate}
\item $S(x+k,y+k)=S(x,y)$ for every $k\in\mathbb Z^d$;
\item for some $\lambda>0$,
\[
 \xi^{\mathsf T}D_{12}S(x,y)\xi\leq-\lambda|\xi|^2
 \qquad(x,y,\xi\in\mathbb R^d);
\]
\item the lift is given by
\begin{equation}\label{eq:unperturbed-generating-equations}
 \widetilde F(x,p)=(y,P)
 \quad\Longleftrightarrow\quad
 \begin{cases}
 p=-D_1S(x,y),\\
 P=D_2S(x,y).
 \end{cases}
\end{equation}
\end{enumerate}
We call $S$ a \emph{generating function} for $\widetilde F$.
\end{definition}

The second condition is the uniform twist condition. The mixed
Hessian $D_{12}S$ need not be symmetric. In dimension one, the condition implies
that the next position is strictly increasing in the initial momentum.
The diagonal periodicity gives exact symplecticity: with
$\theta=\sum_{j=1}^d p_j\,dx_j$, the one-form $F^*\theta-\theta$
is exact.

\begin{lemma}\label{lem:global-map}
Every $C^2$ function $S$ satisfying the diagonal periodicity and uniform
twist conditions defines a global $C^1$ diffeomorphism
$\widetilde F$ by~\eqref{eq:unperturbed-generating-equations}. It satisfies
\[
 \widetilde F(x+k,p)=\widetilde F(x,p)+(k,0),
 \qquad k\in\mathbb Z^d,
\]
and induces an exact symplectic twist map $F$ on
$\mathbb T^d\times\mathbb R^d$.
For each $x\in\mathbb R^d$, the maps
\[
 \mathbb R^d\ni p\longmapsto
 \pi_x\bigl(\widetilde F^{\pm1}(x,p)\bigr)\in\mathbb R^d
\]
are global $C^1$ diffeomorphisms, where $\pi_x(x,p)=x$.
\end{lemma}

For the variational arguments, we impose the following assumptions.

\begin{assumption}\label{ass:generating-function}
The function $S\in C^2(\mathbb R^d\times\mathbb R^d)$ satisfies,
for some $K>\lambda>0$,
\begin{align}
 S(x+k,y+k)&=S(x,y)
 &&(x,y\in\mathbb R^d,\ k\in\mathbb Z^d),
 \label{eq:diagonal-periodicity}\\
 \xi^{\mathsf T}D_{12}S(x,y)\xi&\leq-\lambda|\xi|^2
 &&(x,y,\xi\in\mathbb R^d),
 \label{eq:uniform-twist}\\
 \xi^{\mathsf T}D_{ii}S(x,y)\xi&\leq K|\xi|^2
 &&(x,y,\xi\in\mathbb R^d,\ i=1,2).
 \label{eq:endpoint-bounds}
\end{align}
\end{assumption}

Let $F$ be generated by a function $S$ satisfying
Assumption~\ref{ass:generating-function}, and let $V\in C^\infty(\mathbb T^d)$. The function $S_V(x,y)=S(x,y)-V(x)$
is still diagonally periodic and has the same mixed derivative as $S$.
It satisfies Assumption~\ref{ass:generating-function} after increasing $K$
if necessary. Its map $F_V$ is given by
\eqref{eq:generating-equations}. Comparing the two sets of generating
equations yields
\begin{equation}\label{eq:potential-kick}
 F_V=F\circ K_V,\qquad K_V(x,p)=(x,p-DV(x)).
\end{equation}

For integers $a<b$, the action of a finite lifted configuration is
\[
 \sum_{n=a}^{b-1}\bigl(S(x_n,x_{n+1})-V(x_n)\bigr).
\]
With $x_a$ and $x_b$ fixed, the critical point equations are
\begin{equation}\label{eq:discrete-orbit-equations}
 D_2S(x_{n-1},x_n)+D_1S(x_n,x_{n+1})-DV(x_n)=0,
 \qquad a<n<b.
\end{equation}
Consequently, a complete lifted configuration satisfies these equations
at every integer time if and only if, with
$p_n=D_2S(x_{n-1},x_n)$, the sequence $(x_n,p_n)$ is an orbit of
$\widetilde F_V$. This correspondence will be used to construct orbits
by finite action minimization.

\subsection{The averaged kinetic function}
\label{subsec:periodic-estimates}

Let $S$ satisfy Assumption~\ref{ass:generating-function}.
Define
\begin{equation}\label{eq:kinetic-decomposition}
 L(x,v)=S(x,x+v),\qquad
 h(v)=\int_{\mathbb T^d}L(x,v)\,dx,\qquad
 R(x,v)=L(x,v)-h(v).
\end{equation}
Here $dx$ is the normalized Haar measure on $\mathbb T^d$.

\begin{lemma}\label{lem:averaged-kinetic}
The function $h$ belongs to $C^2$ and satisfies
\begin{equation}\label{eq:averaged-convexity}
 \lambda I\leq D^2h(v)\leq KI.
\end{equation}
The remainder and its derivative with respect to $v$ satisfy
\begin{equation}\label{eq:remainder-bounds}
 \|R\|_\infty\leq\frac d2(K-\lambda),\qquad
 \|D_vR\|_\infty\leq\sqrt d\,(K-\lambda).
\end{equation}
Moreover, $Dh$ is a global $C^1$ diffeomorphism, and its inverse
$(Dh)^{-1}$ satisfies
\begin{equation}\label{eq:inverse-kinetic}
 K^{-1}I\leq D\bigl[(Dh)^{-1}\bigr](c)\leq\lambda^{-1}I.
\end{equation}
Adding a periodic function $J(x)$ to $S$ changes $h$ by a constant
and leaves $Dh$, $(Dh)^{-1}$, $D_vR$, and $D_{22}S$ unchanged.
\end{lemma}

\begin{proof}
Differentiation under the integral gives
\[
 D^2h(v)=\int_{\mathbb T^d}D_{22}S(x,x+v)\,dx.
\]
Since the vector field $x\mapsto D_2S(x,x+v)$ is periodic, its
derivative has zero integral. Thus
\[
 \int_{\mathbb T^d}(D_{21}S+D_{22}S)(x,x+v)\,dx=0
\]
and
\[
 D^2h(v)=-\int_{\mathbb T^d}D_{21}S(x,x+v)\,dx.
\]
The lower bound in \eqref{eq:averaged-convexity} follows from the
twist inequality, and the upper bound follows from $D_{22}S\leq KI$.

For every fixed $v$, differentiating $L(x,v)=S(x,x+v)$ twice
with respect to $x$ gives
\[
 D_x^2L
 =D_{11}S+D_{12}S+D_{21}S+D_{22}S
 \leq2(K-\lambda)I.
\]
For each coordinate vector $e_i$, Taylor's formula gives
\begin{align*}
 L(x+e_i,v)-L(x,v)
 &=\partial_{x_i}L(x,v)
   +\int_0^1(1-t)\partial_{x_i}^2L(x+te_i,v)\,dt\\
 &\leq\partial_{x_i}L(x,v)+(K-\lambda),\\
 L(x-e_i,v)-L(x,v)
 &=-\partial_{x_i}L(x,v)
   +\int_0^1(1-t)\partial_{x_i}^2L(x-te_i,v)\,dt\\
 &\leq-\partial_{x_i}L(x,v)+(K-\lambda).
\end{align*}
By periodicity, both left-hand sides are zero. The first inequality
gives $\partial_{x_i}L(x,v)\geq-(K-\lambda)$, and the second gives
$\partial_{x_i}L(x,v)\leq K-\lambda$. Thus
\[
 |\partial_{x_i}L(x,v)|\leq K-\lambda
 \qquad(1\leq i\leq d).
\]
Taking the Euclidean norm gives
$|D_xL|\leq\sqrt d\,(K-\lambda)$.
Since the torus has diameter $\sqrt d/2$, comparison with the
$x$-average gives the first bound in \eqref{eq:remainder-bounds}.

For the second bound, fix $v$ and set
\[
 G_v(x)=D_2S(x,x+v).
\]
This periodic vector field is $C^1$, with
$DG_v=D_{21}S+D_{22}S$. By Assumption~\ref{ass:generating-function},
\[
 \xi^{\mathsf T}DG_v(x)\xi\leq(K-\lambda)|\xi|^2.
\]
For $x,t\in\mathbb R^d$, integration and periodicity give
\[
 \bigl(G_v(x)-G_v(t)\bigr)\cdot z\leq(K-\lambda)|z|^2
 \qquad\bigl(z=x-t+k,\ k\in\mathbb Z^d\bigr).
\]
Since the flat torus has diameter $\sqrt d/2$, we can choose
$k\in\mathbb Z^d$ so that
\[
 \left|z-\frac{G_v(x)-G_v(t)}{2(K-\lambda)}\right|
 \leq\frac{\sqrt d}{2}.
\]
Completing the square gives
\begin{align*}
 0&\leq(K-\lambda)|z|^2-\bigl(G_v(x)-G_v(t)\bigr)\cdot z\\
  &=(K-\lambda)
    \left|z-\frac{G_v(x)-G_v(t)}{2(K-\lambda)}\right|^2
    -\frac{|G_v(x)-G_v(t)|^2}{4(K-\lambda)}\\
  &\leq\frac d4(K-\lambda)
    -\frac{|G_v(x)-G_v(t)|^2}{4(K-\lambda)}.
\end{align*}
Therefore,
\[
 |G_v(x)-G_v(t)|\leq\sqrt d\,(K-\lambda).
\]
Averaging in $t$ now yields
\[
 |D_vR(x,v)|
 =\left|G_v(x)-\int_{\mathbb T^d}G_v(t)\,dt\right|
 \leq\sqrt d\,(K-\lambda).
\]

By~\eqref{eq:averaged-convexity}, for every $u,v\in\mathbb R^d$,
\begin{align*}
 \bigl(Dh(v)-Dh(u)\bigr)\cdot(v-u)
 &=\int_0^1(v-u)^{\mathsf T}D^2h\bigl(u+s(v-u)\bigr)(v-u)\,ds\\
 &\geq\lambda|v-u|^2.
\end{align*}
This inequality is the strong monotonicity of $Dh$.
It implies that $Dh$ is injective and that
$|Dh(v)|\geq\lambda|v|-|Dh(0)|$. Thus $Dh$ is proper.
Since $D^2h$ is invertible everywhere, $Dh$ is locally invertible.
Its image is therefore open and closed in $\mathbb R^d$, so $Dh$
is a global $C^1$ diffeomorphism.
Inverting the positive definite Hessian gives \eqref{eq:inverse-kinetic}.
Finally, adding $J(x)$ changes $h$ by $\int_{\mathbb T^d}J$ and
$R$ by $J-\int_{\mathbb T^d}J$. This proves the remaining assertion.
\end{proof}

\subsection{Minimizing configurations under an affine force}
\label{subsec:affine-minimizers}

Let $S$ satisfy Assumption~\ref{ass:generating-function}.

For $J\in C^\infty(\mathbb T^d)$ and $p^*\in\mathbb R^d$, consider
the auxiliary generating function
\begin{equation}\label{eq:affine-generator}
 \mathcal{S}_{J,p^*}(x,y)=S(x,y)+J(x)+p^*\cdot x.
\end{equation}
For every $k\in\mathbb Z^d$,
\[
 \mathcal{S}_{J,p^*}(x+k,y+k)
 =\mathcal{S}_{J,p^*}(x,y)+p^*\cdot k.
\]
Hence $\mathcal{S}_{J,p^*}$ is diagonally periodic if and only if $p^*=0$.
Its first derivatives are diagonally periodic for every $p^*$.
The generating equations are
\begin{equation}\label{eq:affine-map}
 p=-D_1S(x,y)-DJ(x)-p^*,\qquad P=D_2S(x,y).
\end{equation}
These equations define the global $C^1$ diffeomorphism
\[
 \mathcal{F}_{J,p^*}(x,p)=F\bigl(x,p+DJ(x)+p^*\bigr)
\]
of $\mathbb T^d\times\mathbb R^d$, where $F$ is generated by $S$.
Its lift to $\mathbb R^d\times\mathbb R^d$ is
\[
 \widetilde{\mathcal{F}_{J,p^*}}(x,p)
 =\widetilde F\bigl(x,p+DJ(x)+p^*\bigr).
\]
It is generated by $\mathcal{S}_{J,p^*}$ and has the same uniform
twist condition, since $D_{12}\mathcal{S}_{J,p^*}=D_{12}S$.
The momentum translation above changes $\theta$ by the closed form
\[
 dJ+\sum_{j=1}^d p_j^*\,dx_j.
\]
Its integral along the $j$th coordinate circle is
\[
 \int_0^1\bigl(\partial_{x_j}J(x+te_j)+p_j^*\bigr)\,dt
 =J(x+e_j)-J(x)+p_j^*=p_j^*.
\]
Since $F$ is exact symplectic, composition with $F$ adds an exact
form and leaves these integrals unchanged.
Thus $\mathcal{F}_{J,p^*}$ is symplectic. It is exact symplectic if and
only if $p^*=0$.
For $p^*\ne0$, $\mathcal{F}_{J,p^*}$ is therefore not an exact symplectic twist map
in the sense of Definition~\ref{def:twist-map}.

\begin{definition}[Globally minimizing configurations]
\label{def:globally-minimizing}
A sequence $(x_n)_{n\in\mathbb Z}\subset\mathbb R^d$ is
\emph{globally minimizing for $\mathcal{S}_{J,p^*}$} if, for every
pair of integers $a<b$,
\[
 \sum_{n=a}^{b-1}\mathcal{S}_{J,p^*}(x_n,x_{n+1})
 =\min_{\substack{y_a,\ldots,y_b\in\mathbb R^d\\
                   y_a=x_a,\;y_b=x_b}}
   \sum_{n=a}^{b-1}\mathcal{S}_{J,p^*}(y_n,y_{n+1}).
\]
\end{definition}

\begin{lemma}\label{lem:affine-minimizer}
For every smooth periodic $J$ and every $p^*\in\mathbb R^d$, the lift
$\widetilde{\mathcal{F}_{J,p^*}}$ has an orbit $(x_n,p_n)_{n\in\mathbb Z}$.
The configuration $(x_n)_{n\in\mathbb Z}$ is globally minimizing for
$\mathcal{S}_{J,p^*}$, and
the momenta satisfy
\begin{equation}\label{eq:affine-tube}
 |p_n-np^*|\le \frac{2\sqrt d\,K^2}{\lambda}\qquad(n\in\mathbb Z).
\end{equation}
This bound is independent of both $J$ and $p^*$.
\end{lemma}

\begin{proof}
\noindent\emph{Step 1. Choosing the endpoints and the finite action.}
Fix a positive integer $N$ and $c_0\in\mathbb R^d$, and set
$c_j=c_0+jp^*$ for $1\leq j\leq N$.
The vector $c_j$ is the target value of $p_j$.
By Lemma~\ref{lem:averaged-kinetic}, the function
\[
 v\longmapsto h(v)-c_j\cdot v
\]
has the unique minimizer $(Dh)^{-1}(c_j)$.
We therefore choose
\[
 x_0=0,\qquad x_N=\sum_{j=1}^{N}(Dh)^{-1}(c_j).
\]
With $x_0$ and $x_N$ fixed, minimize
\begin{equation}\label{eq:finite-cohomology-action}
 \mathcal B_N(x_0,\ldots,x_N)
 =\sum_{j=1}^{N}
 \left[S(x_{j-1},x_j)+J(x_{j-1})-c_j\cdot(x_j-x_{j-1})\right].
\end{equation}
The decomposition~\eqref{eq:kinetic-decomposition} writes each summand
as $h(v_j)-c_j\cdot v_j+R(x_{j-1},v_j)+J(x_{j-1})$, where
$v_j=x_j-x_{j-1}$. The functions $R$ and $J$ are bounded, and $h$
has a uniform quadratic lower bound. Thus $\mathcal B_N$ is coercive
and attains a global minimum.

Since $c_{j+1}-c_j=p^*$, summation by parts gives
\[
 \mathcal B_N
 =\sum_{j=1}^{N}\mathcal{S}_{J,p^*}(x_{j-1},x_j)
  +c_0\cdot x_0-c_N\cdot x_N.
\]
The last two terms depend only on the endpoints.
Hence every minimizer of $\mathcal B_N$, and each of its subsegments,
minimizes the $\mathcal{S}_{J,p^*}$ action with fixed lifted endpoints.
For such a minimizer, set
\[
 p_j=D_2S(x_{j-1},x_j)\qquad(1\leq j\leq N).
\]
The interior Euler--Lagrange equations are
\[
 D_2S(x_{j-1},x_j)+D_1S(x_j,x_{j+1})+DJ(x_j)=-p^*
 \qquad(1\leq j<N).
\]
These are the orbit equations of $\widetilde{\mathcal{F}_{J,p^*}}$.

\smallskip
\noindent\emph{Step 2. Bounding the differences between momentum deviations.}
Write
\[
 d_j=p_j-c_j\qquad(1\leq j\leq N).
\]
For $1\leq i<\ell\leq N$ and $k\in\mathbb Z^d$, translate the block
$x_i,\ldots,x_{\ell-1}$ by $k$. Define the change in action by
\[
 \begin{aligned}
 \Delta\mathcal B_N(k)
 ={}&\mathcal B_N(x_0,\ldots,x_{i-1},x_i+k,\ldots,x_{\ell-1}+k,
                   x_\ell,\ldots,x_N)\\
 &-\mathcal B_N(x_0,\ldots,x_N).
 \end{aligned}
\]
All summands except those with indices $i$ and $\ell$ remain
unchanged. By diagonal periodicity,
$S(x_{\ell-1}+k,x_\ell)=S(x_{\ell-1},x_\ell-k)$.
We apply Taylor's formula with integral remainder to these two terms
in their second variables. The bound $D_{22}S\leq KI$ gives an upper
bound $K|k|^2/2$ for each remainder.
Minimality gives
\[
 0\leq\Delta\mathcal B_N(k)
 \leq(d_i-d_\ell)\cdot k+K|k|^2.
\]
Taking $k=e_h$ and $k=-e_h$ for each coordinate vector $e_h$ gives
\begin{equation}\label{eq:deviation-oscillation}
 |d_i-d_j|\leq\sqrt d\,K\qquad(1\leq i,j\leq N).
\end{equation}
Thus all deviations stay within a fixed distance of $d_1$.
It remains to bound $d_1$.

\smallskip
\noindent\emph{Step 3. Using the total displacement to bound $d_1$.}
Since $p_j=c_j+d_j$, \eqref{eq:kinetic-decomposition} gives
\[
 Dh(v_j)=c_j+d_j-D_vR(x_{j-1},v_j).
\]
We now use the endpoint choice from Step~1. It gives the identity
\[
 \sum_{j=1}^{N}v_j
 =\sum_{j=1}^{N}(Dh)^{-1}(c_j).
\]
The bounds~\eqref{eq:inverse-kinetic},
\eqref{eq:deviation-oscillation}, and~\eqref{eq:remainder-bounds} yield
\begin{align*}
 \frac{|d_1|^2}{K}
 &\leq\frac1N\sum_{j=1}^{N}
 d_1\cdot\bigl[(Dh)^{-1}(c_j+d_1)-(Dh)^{-1}(c_j)\bigr]\\
 &=\frac1N\sum_{j=1}^{N}
 d_1\cdot\bigl[(Dh)^{-1}(c_j+d_1)-v_j\bigr]\\
 &\leq\frac{|d_1|}{N\lambda}\sum_{j=1}^{N}
 \bigl|d_1-d_j+D_vR(x_{j-1},v_j)\bigr|\\
 &\leq\frac{|d_1|}{N\lambda}\sum_{j=1}^{N}
 \bigl(|d_j-d_1|+|D_vR(x_{j-1},v_j)|\bigr)\\
 &\leq\frac{\sqrt d\,(2K-\lambda)}{\lambda}|d_1|.
\end{align*}
The equality uses the prescribed total displacement.
It follows that $|d_1|\leq\sqrt d\,K(2K-\lambda)/\lambda$.
Combining this bound with~\eqref{eq:deviation-oscillation}, we obtain
\begin{equation}\label{eq:finite-affine-tube}
 |p_j-c_j|
 \leq\sqrt d\,K+\frac{\sqrt d\,K(2K-\lambda)}{\lambda}
 =\frac{2\sqrt d\,K^2}{\lambda}\qquad(1\leq j\leq N).
\end{equation}
The two estimates have different roles:
integer translations control the differences $d_i-d_j$, and the
endpoint constraint controls $d_1$.
Together they give a bound independent of $N$, $c_0$, $J$, and $p^*$.

\smallskip
\noindent\emph{Step 4. Passing to an orbit.}
We use a compactness argument similar to that in
\cite[Section~4.3, Step~5 and Lemma~4.5]{GambaudoGuiraudPetite2006}.
Apply the finite construction on the interval $-N,\ldots,N$.
Set $c_j=jp^*$ for $-N+1\leq j\leq N$, fix $x_{-N}=0$, and choose
\[
 x_N=\sum_{j=-N+1}^{N}(Dh)^{-1}(c_j).
\]
The estimate~\eqref{eq:finite-affine-tube} becomes
\[
 |p_n-np^*|\leq\frac{2\sqrt d\,K^2}{\lambda}
 \qquad(-N+1\leq n\leq N).
\]
In particular, $p_0$ stays in a fixed ball about $0$.
Translate each configuration by an integer vector so that $x_0$ lies
in a fixed closed fundamental cube.
This translation changes the action of every competitor by the same
constant and preserves minimality.
Passing to a subsequence, the phase points $(x_0,p_0)$ at time zero converge.
Since $\widetilde{\mathcal{F}_{J,p^*}}$ and its inverse are continuous,
every fixed positive or negative iterate also converges.
The limit is an orbit of $\widetilde{\mathcal{F}_{J,p^*}}$ satisfying
\eqref{eq:affine-tube} for every $n\in\mathbb Z$.

The limit configuration is globally minimizing.
Indeed, suppose that one of its finite segments has a competitor with
strictly smaller action. Keep the interior points of this competitor
and replace its endpoints by those of an approximating segment.
By continuity, the action remains smaller for all sufficiently close
approximations. This contradicts minimality.
\end{proof}

\subsection{Small potentials and the removal of an affine force}
\label{subsec:potential-tools}

We first construct a small periodic potential and a periodic correction
to an affine function. We then use finite differences to remove the
affine force along minimizing configurations.

\begin{figure}[htbp]
\centering
\input{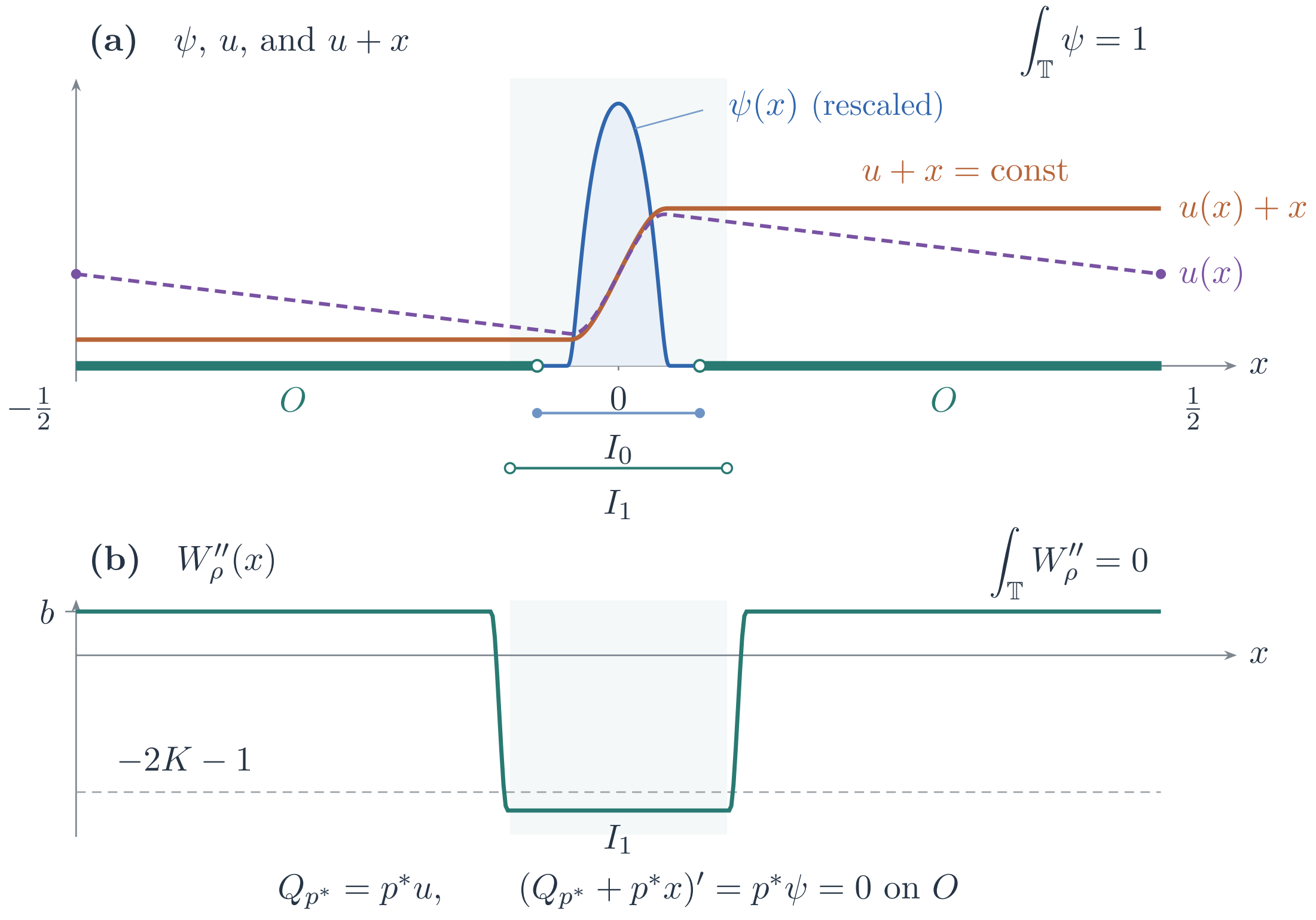}
\caption{The one-dimensional construction in Lemma~\ref{lem:potential-correction}.
The blue curve is rescaled vertically.
The endpoints $-1/2$ and $1/2$ are identified.}
\label{fig:small-potential-correction}
\end{figure}

\FloatBarrier

\begin{lemma}
\label{lem:potential-correction}
Let $d\geq1$, $K>0$, $\varepsilon>0$, and $0<\alpha<1$.
There exist $W\in C^\infty(\mathbb T^d)$,
$u\in C^\infty(\mathbb T)$, an open interval $I_1\subset\mathbb T$,
a closed interval $I_0\subset I_1$, and $0<r\leq1$ such that
\begin{equation}\label{eq:negative-strip-potential}
 \begin{gathered}
 \|W\|_{C^{1,\alpha}}<\varepsilon,\\
 \partial_{x_j}^2W(x)\leq-2K-1
 \quad\bigl(x\in\mathbb T^{j-1}\times I_1\times\mathbb T^{d-j},\ 1\leq j\leq d\bigr).
 \end{gathered}
\end{equation}
Set $O=(\mathbb T\setminus I_0)^d$. For every $p^*\in\mathbb R^d$
with $|p^*|\leq r$, the function
\[
 Q_{p^*}(x)=\sum_{j=1}^d p_j^* u(x_j)
\]
satisfies
\begin{equation}\label{eq:auxiliary-correction}
 DQ_{p^*}(x)=-p^*\quad(x\in O),\qquad
 \|D^2Q_{p^*}\|_\infty\leq\frac14.
\end{equation}
\end{lemma}

\begin{proof}
Fix a nonnegative function $\varphi\in C_c^\infty((-1/4,1/4))$
which is at most one and equals one on $[-1/8,1/8]$.
For small $\rho>0$, use a coordinate chart about $0\in\mathbb T$
and extend $\varphi(s/\rho)$ by zero to the circle. Set
\[
 b=2(K+1)\rho\int_{\mathbb R}\varphi(s)\,ds.
\]
Choose the mean-zero smooth periodic function $w_\rho$ satisfying
\begin{equation}\label{eq:spike-second-derivative}
 w_\rho''(s)=b-2(K+1)\varphi(s/\rho).
\end{equation}
The right-hand side has mean zero, so such a function exists.
Choose $\rho$ small enough that $b\leq1$.
Since $\|w_\rho''\|_{L^1}\leq2b$ and both $w_\rho$ and $w_\rho'$
have mean zero, we have
\[
 \|w_\rho'\|_\infty\leq2b,\qquad
 \|w_\rho\|_\infty\leq2b.
\]
Define
\begin{equation}\label{eq:final-spike-potential}
 W_\rho(x)=\sum_{j=1}^d w_\rho(x_j).
\end{equation}
Then
\begin{equation}\label{eq:spike-norms}
 \|W_\rho\|_\infty\leq2db,\qquad
 \|DW_\rho\|_\infty\leq2\sqrt d\,b,\qquad
 \|D^2W_\rho\|_\infty\leq2K+3.
\end{equation}
For $x,y\in\mathbb T^d$,
\[
 |DW_\rho(x)-DW_\rho(y)|
 \leq\min\bigl\{2\|DW_\rho\|_\infty,
       \|D^2W_\rho\|_\infty\operatorname{dist}_{\mathbb T^d}(x,y)\bigr\}.
\]
It follows that
\begin{equation}\label{eq:spike-holder}
 [DW_\rho]_\alpha
 \leq(2\|DW_\rho\|_\infty)^{1-\alpha}
       \|D^2W_\rho\|_\infty^\alpha
 \leq(4\sqrt d\,b)^{1-\alpha}(2K+3)^\alpha.
\end{equation}
Thus $\|W_\rho\|_{C^{1,\alpha}}\to0$ as $\rho\to0$.
Choose $\rho$ so that this norm is less than $\varepsilon$, and put
$W=W_\rho$ and $I_1=\{s\in\mathbb T:|s|<\rho/8\}$.
On $I_1$, we have $w_\rho''=b-2(K+1)\leq-2K-1$.
This proves~\eqref{eq:negative-strip-potential}.

Choose a closed interval $I_0\subset I_1$ and a nonnegative function
$\psi\in C^\infty(\mathbb T)$ such that
\[
 \int_{\mathbb T}\psi(s)\,ds=1,\qquad
 \operatorname{supp}\psi\subset I_0.
\]
Choose
\[
 u\in C^\infty(\mathbb T),\qquad u'(s)=\psi(s)-1
 \quad(s\in\mathbb T).
\]
Then $u'=-1$ on $\mathbb T\setminus I_0$.
Choose $0<r\leq1$ such that $r\|u''\|_\infty\leq1/4$.
For $|p^*|\leq r$, the definition of $Q_{p^*}$ gives
$DQ_{p^*}=-p^*$ on $O$, and
\[
 \|D^2Q_{p^*}\|_\infty
 \leq |p^*|\|u''\|_\infty\leq\frac14.
\]
\end{proof}

\begin{lemma}
\label{lem:exclude-and-cancel}
Let $S$ satisfy Assumption~\ref{ass:generating-function}.
Let $W\in C^\infty(\mathbb T^d)$, let $I_0\subset\mathbb T$ be a closed
interval, and let $r>0$. Suppose that
\[
 \partial_{x_j}^2W(x)\leq-2K-1
 \quad\bigl(x\in\mathbb T^{j-1}\times I_0\times\mathbb T^{d-j},\ 1\leq j\leq d\bigr).
\]
Set $O=(\mathbb T\setminus I_0)^d$. Suppose that, for every
$p^*\in\mathbb R^d$ with $|p^*|\leq r$, there is a function
$Q_{p^*}\in C^\infty(\mathbb T^d)$ such that
\[
 DQ_{p^*}(x)=-p^*\quad(x\in O),\qquad
 \|D^2Q_{p^*}\|_\infty\leq\frac14.
\]
Then, for every $|p^*|\leq r$, there is a point
$z_{p^*}\in\mathbb T^d\times\mathbb R^d$ such that
\begin{equation}\label{eq:local-drift}
 \sup_{n\in\mathbb Z}
 \big|\pi_pF_{-W}^{\,n}(z_{p^*})-np^*\big|\leq \frac{2\sqrt d\,K^2}{\lambda}.
\end{equation}
\end{lemma}

\begin{proof}
Fix $|p^*|\leq r$ and choose $Q_{p^*}$ as in the statement.
Apply Lemma~\ref{lem:affine-minimizer} with $J=W+Q_{p^*}$.
It gives an orbit $(x_n,p_n)_{n\in\mathbb Z}$ of
$\widetilde{\mathcal{F}_{W+Q_{p^*},p^*}}$ whose
configuration globally minimizes
\[
 \mathcal{S}_{W+Q_{p^*},p^*}(x,y)
 =S(x,y)+W(x)+Q_{p^*}(x)+p^*\cdot x
\]
and whose momenta satisfy~\eqref{eq:affine-tube}.

Fix $n\in\mathbb Z$. Keep $x_{n-1},x_{n+1}$ fixed and vary $x_n$.
The part of the action that depends on this point is
\[
 \Phi(x)=S(x_{n-1},x)+S(x,x_{n+1})+W(x)+Q_{p^*}(x)+p^*\cdot x.
\]
The bounds on $D_{11}S$, $D_{22}S$, $\partial_{x_j}^2W$, and $D^2Q_{p^*}$ give
\begin{equation}\label{eq:uniform-curvature-gap}
 \partial_{x_j}^2\Phi(x)
 \leq2K-(2K+1)+\frac14=-\frac34
 \qquad(x_j\bmod1\in I_0).
\end{equation}
If $(x_n)_j\bmod1\in I_0$ for some $j$, then
$\partial_{x_j}^2\Phi(x_n)<0$. Minimality gives
$\partial_{x_j}^2\Phi(x_n)\geq0$, a contradiction.
Thus $x_n\bmod\mathbb Z^d$ lies in $O$ for every $n\in\mathbb Z$.

On $O$, we have $DQ_{p^*}=-p^*$. Hence the momentum equations become
\[
 \begin{split}
 p_n&=-D_1S(x_n,x_{n+1})-DW(x_n)-DQ_{p^*}(x_n)-p^*\\
    &=-D_1S(x_n,x_{n+1})-DW(x_n),\\
 p_{n+1}&=D_2S(x_n,x_{n+1}).
 \end{split}
\]
These are the generating equations for $S+W$.
The projected phase sequence is therefore an orbit of $F_{-W}$.
Set
\[
 z_{p^*}=(x_0\bmod\mathbb Z^d,p_0)\in\mathbb T^d\times\mathbb R^d.
\]
Then \eqref{eq:affine-tube} gives~\eqref{eq:local-drift}.
\end{proof}

\begin{lemma}
\label{lem:drift-persistence}
Let $S$ satisfy Assumption~\ref{ass:generating-function}.
Let $W\in C^\infty(\mathbb T^d)$, let $I_0\subset I_1\subset\mathbb T$
be intervals with $I_0$ closed and $I_1$ open, and let $r>0$.
Suppose that
\[
 \partial_{x_j}^2W(x)\leq-2K-1
 \quad\bigl(x\in\mathbb T^{j-1}\times I_1\times\mathbb T^{d-j},\ 1\leq j\leq d\bigr).
\]
Set $O=(\mathbb T\setminus I_0)^d$. Suppose that, for every
$p^*\in\mathbb R^d$ with $|p^*|\leq r$, there is a function
$Q_{p^*}\in C^\infty(\mathbb T^d)$ such that
\[
 DQ_{p^*}(x)=-p^*\quad(x\in O),\qquad
 \|D^2Q_{p^*}\|_\infty\leq\frac14.
\]
Choose
\begin{equation}\label{eq:persistence-radius}
 0<t<\operatorname{dist}(I_0,\partial I_1),\qquad
 \delta=\frac{t^2}{16}.
\end{equation}
For every $E\in C^\infty(\mathbb T^d)$ with $\|E\|_\infty<\delta$
and every $p^*\in\mathbb R^d$ with $|p^*|\leq r$, there is a point
$z_{E,p^*}\in\mathbb T^d\times\mathbb R^d$ such that
\begin{equation}\label{eq:local-persistent-drift}
 \sup_{n\in\mathbb Z}
 \big|\pi_pF_{-W+E}^{\,n}(z_{E,p^*})-np^*\big|
 \leq\frac{2\sqrt d\,K^2}{\lambda}.
\end{equation}
The radius $r$ and the bound are independent of $E$.
\end{lemma}

\begin{proof}
Fix $E\in C^\infty(\mathbb T^d)$ with $\|E\|_\infty<\delta$ and
$|p^*|\leq r$, and choose $Q_{p^*}$ as in the statement.
Apply Lemma~\ref{lem:affine-minimizer} with $J=W-E+Q_{p^*}$.
It gives an orbit $(x_n,p_n)_{n\in\mathbb Z}$ whose configuration globally minimizes
\begin{equation}\label{eq:perturbed-affine-action}
 \mathcal{S}_{W-E+Q_{p^*},p^*}(x,y)
 =S(x,y)+W(x)-E(x)+Q_{p^*}(x)+p^*\cdot x
\end{equation}
and whose momenta satisfy~\eqref{eq:affine-tube}.

Fix $n\in\mathbb Z$ and suppose that $(x_n)_j\bmod1\in I_0$ for some $j$.
Keep $x_{n-1},x_{n+1}$ fixed and vary $x_n$.
The part of the action that depends on this point is $\Phi(x)-E(x)$, where
\[
 \Phi(x)=S(x_{n-1},x)+S(x,x_{n+1})+W(x)+Q_{p^*}(x)+p^*\cdot x.
\]
For $|s|\leq t$, the $j$th coordinate of $x_n+se_j$, modulo one, remains in $I_1$,
so the bounds on $D_{11}S$, $D_{22}S$, $\partial_{x_j}^2W$, and $D^2Q_{p^*}$ give
\[
 \partial_{x_j}^2\Phi(x_n+se_j)
 \leq2K-(2K+1)+\frac14=-\frac34.
\]
Replace only $x_n$ by $x_n+te_j$ or $x_n-te_j$ and keep all other
points fixed. Both variations keep the endpoints $x_{n-1}$ and
$x_{n+1}$ fixed. Minimality gives
\[
 \begin{aligned}
  (\Phi-E)(x_n+te_j)&\geq(\Phi-E)(x_n),\\
  (\Phi-E)(x_n-te_j)&\geq(\Phi-E)(x_n).
 \end{aligned}
\]
We show that these two inequalities cannot both hold.

Taylor's formula with integral remainder gives
\[
 \Phi(x_n\pm te_j)-\Phi(x_n)
 =\pm t\,\partial_{x_j}\Phi(x_n)
  +\int_0^t(t-s)\,\partial_{x_j}^2\Phi(x_n\pm se_j)\,ds.
\]
Adding the two identities cancels the first derivative terms. Hence
\begin{equation}\label{eq:finite-difference-integral}
 \begin{aligned}
  &\Phi(x_n+te_j)+\Phi(x_n-te_j)-2\Phi(x_n)\\
  &\qquad=\int_{-t}^{t}(t-|s|)\,\partial_{x_j}^2\Phi(x_n+se_j)\,ds\\
  &\qquad\leq-\frac34\int_{-t}^{t}(t-|s|)\,ds
  =-\frac{3t^2}{4}.
 \end{aligned}
\end{equation}
The contribution of $E$ satisfies
\[
 \big|E(x_n+te_j)+E(x_n-te_j)-2E(x_n)\big|
 \leq4\|E\|_\infty.
\]
Adding the two inequalities from minimality and using
$\|E\|_\infty<\delta=t^2/16$, we obtain
\[
 \begin{aligned}
  0&\leq(\Phi-E)(x_n+te_j)+(\Phi-E)(x_n-te_j)-2(\Phi-E)(x_n)\\
   &\leq-\frac{3t^2}{4}+4\|E\|_\infty
   <-\frac{t^2}{2}<0.
 \end{aligned}
\]
This contradiction excludes $(x_n)_j\bmod1\in I_0$ for every $j$.
Thus $x_n\bmod\mathbb Z^d\in O$ for every $n\in\mathbb Z$.

The identity $DQ_{p^*}=-p^*$ on $O$ reduces the momentum equations to
\[
 \begin{split}
 p_n&=-D_1S(x_n,x_{n+1})+D(-W+E)(x_n),\\
 p_{n+1}&=D_2S(x_n,x_{n+1}).
 \end{split}
\]
These are the generating equations for $S-(-W+E)$.
The projected phase sequence is therefore an orbit of $F_{-W+E}$.
Set
\[
 z_{E,p^*}=(x_0\bmod\mathbb Z^d,p_0)\in\mathbb T^d\times\mathbb R^d.
\]
Then \eqref{eq:affine-tube} gives~\eqref{eq:local-persistent-drift}.
\end{proof}

\subsection{Repeated approximation and localized generating functions}
\label{subsec:localized-generators}

\begin{lemma}\label{lem:repeated-potentials}
Let $d\geq1$ and $\varepsilon>0$. There is a sequence
$(U_k)_{k\geq1}$ of real trigonometric polynomials on $\mathbb T^d$
such that
\[
 \|U_k\|_{C^2}<2\varepsilon.
\]
For every pair of positive integers $j,N$, there is an integer
$k\geq N$ such that $U_k=U_j$. Moreover, for every
$V\in C^\infty(\mathbb T^d)$ with $\|V\|_{C^2}\leq\varepsilon$,
every integer $r\geq2$, every $\delta>0$, and every $N\geq1$,
there is $k\geq N$ such that
\[
 \|U_k-V\|_{C^r}<\delta.
\]
\end{lemma}

\begin{proof}
List all real trigonometric polynomials with rational coefficients
whose $C^2$ norm is less than $2\varepsilon$, and repeat each
polynomial infinitely often. Approximate $V$ in $C^r$ by a
trigonometric polynomial, then approximate its coefficients by rational
numbers. Since $r\geq2$ and $\|V\|_{C^2}\leq\varepsilon$, sufficiently
close approximations have $C^2$ norm less than $2\varepsilon$.
They therefore belong to the list and occur at arbitrarily large
indices.
\end{proof}

\begin{lemma}\label{lem:localized-generator}
Let $(U_k)$ be smooth periodic functions with
$\sup_k\|U_k\|_{C^2}<\infty$. Let
$\chi_k\in C_c^\infty(\mathbb R^d,[0,1])$ satisfy
$\sup_k\|\chi_k\|_{C^2}<\infty$.
Assume that their supports are pairwise disjoint, that is,
\[
 \operatorname{supp}\chi_i\cap\operatorname{supp}\chi_j=\varnothing
 \qquad(i\ne j).
\]
Assume also that the supports are locally finite: for every compact
set $B\subset\mathbb R^d$,
\[
 \#\{k\geq1:\operatorname{supp}\chi_k\cap B\ne\varnothing\}<\infty.
\]
Then
\begin{equation}\label{eq:localized-generator}
 S(x,y)=\frac12|y-x|^2+\sum_{k\geq1}\chi_k(y-x)U_k(x)
\end{equation}
is smooth, diagonally periodic, and has bounded Hessian. Moreover,
\begin{equation}\label{eq:localized-twist-bound}
 \|D_{12}S+I\|_\infty
 \leq\sup_k\left(
 \|DU_k\|_\infty\|D\chi_k\|_\infty
 +\|U_k\|_\infty\|D^2\chi_k\|_\infty\right).
\end{equation}
If the right-hand side is at most $1/2$, then
\begin{equation}\label{eq:shell-twist-bound}
 \frac12|\xi|^2\leq-\xi^{\mathsf T}D_{12}S(x,y)\xi
 \leq\frac32|\xi|^2.
\end{equation}

In this case, let $V\in C^\infty(\mathbb T^d)$ and put $v=y-x$.
On the support of one term, write $U=U_k$ and $\chi=\chi_k$.
The generating equations of $S-V(x)$ are
\begin{equation}\label{eq:catalogue-generating-equations}
 P-p=\chi(v)DU(x)-DV(x),
 \qquad P=v+U(x)D\chi(v).
\end{equation}
Outside the supports, the equations reduce to
\[
 p=v+DV(x),\qquad P=v.
\]
\end{lemma}

\begin{proof}
Local finiteness gives smoothness, and periodicity of the $U_k$
gives diagonal periodicity. At each point, at most one term contributes.
Direct differentiation on its support gives
\begin{align}
 D_{11}S&=I+\chi D^2U-DU\otimes D\chi
              -D\chi\otimes DU+U D^2\chi,
              \label{eq:shell-D11}\\
 D_{12}S&=-I+DU\otimes D\chi-U D^2\chi,
              \label{eq:shell-D12}\\
 D_{22}S&=I+U D^2\chi.
              \label{eq:shell-D22}
\end{align}
These formulas prove the Hessian bound and
\eqref{eq:localized-twist-bound}. The bound
\eqref{eq:shell-twist-bound} follows by evaluating $D_{12}S+I$
on $\xi$. The first derivatives give
\eqref{eq:catalogue-generating-equations}.
\end{proof}

\Needspace{9cm}
\subsection{Slow passage through annuli}\label{subsec:slow-passage}
\leavevmode\par\nopagebreak

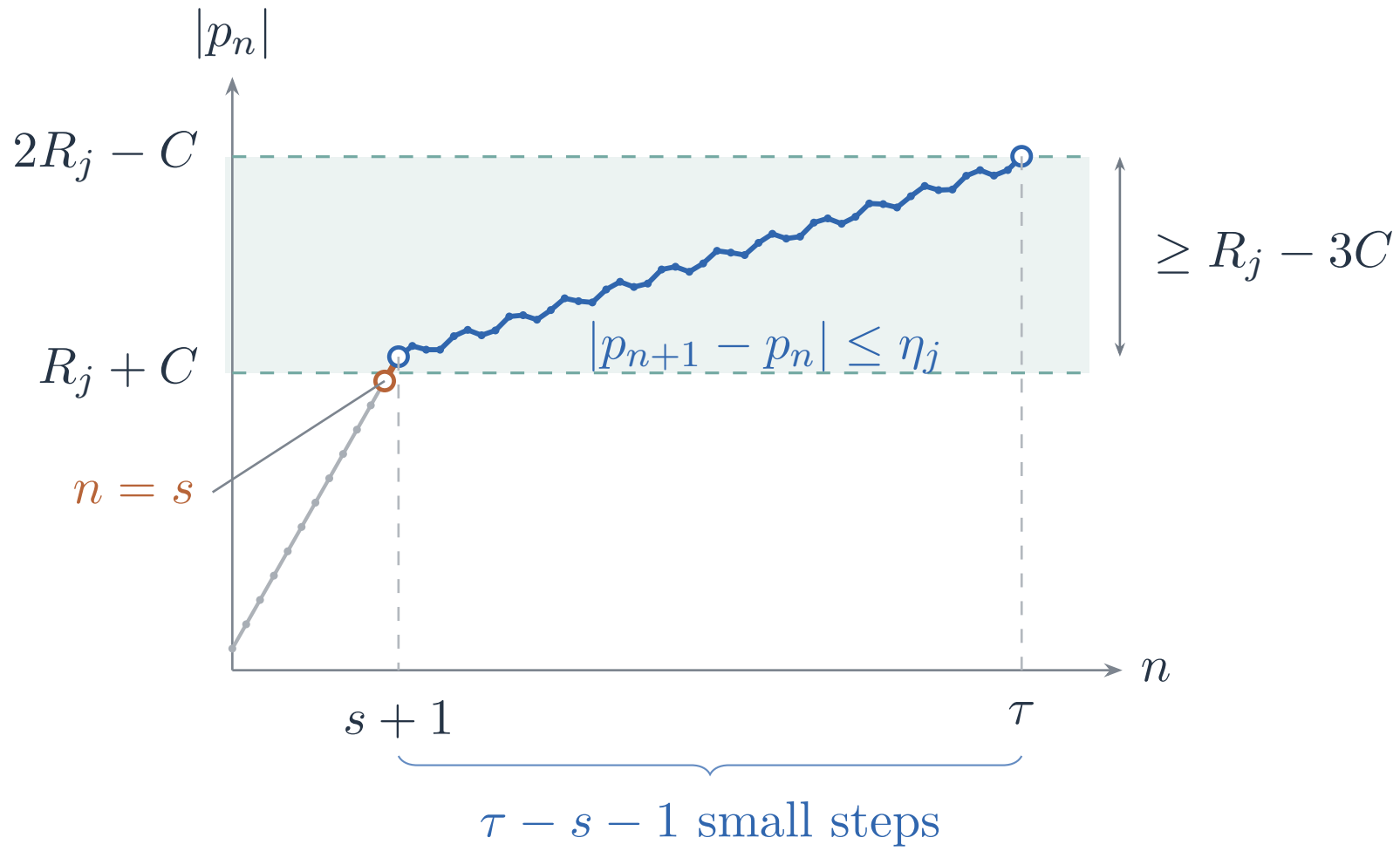
\begin{figure}[htbp]
\centering
\begingroup
\definecolor{ink}{HTML}{263445}
\definecolor{slowcolor}{HTML}{3066AE}
\definecolor{bandcolor}{HTML}{297A72}
\definecolor{entrycolor}{HTML}{B76438}
\begin{tikzpicture}[
 every node/.style={font=\normalsize,text=ink},
 axis/.style={draw=ink!60,line width=.5pt,-{Stealth[length=4pt]}},
 leader/.style={draw=ink!60,line width=.5pt},
 slow/.style={draw=slowcolor,line width=1.1pt,line join=round},
 pre/.style={draw=ink!38,line width=.8pt},
 event/.style={fill=white,draw=slowcolor,line width=.9pt},
 dim/.style={draw=ink!65,line width=.5pt,{Stealth[length=3pt]}-{Stealth[length=3pt]}}]
\fill[bandcolor!9] (9.05000,2.65500) rectangle (15.60000,4.29500);
\draw[axis] (9.10000,0.40000)--(15.85000,0.40000) node[right] {$n$};
\draw[axis] (9.10000,0.40000)--(9.10000,4.90000) node[above] {$|p_n|$};
\draw[bandcolor!65,dashed,line width=.6pt] (9.10000,2.65500)--(15.60000,2.65500);
\node[anchor=east] at (8.98000,2.65500) {$R_j+C$};
\draw[bandcolor!65,dashed,line width=.6pt] (9.10000,4.29500)--(15.60000,4.29500);
\node[anchor=east] at (8.98000,4.29500) {$2R_j-C$};
\draw[pre] (9.10000,0.56400)--(9.20500,0.74850)--(9.31000,0.93300)--(9.41500,1.11750)--(9.52000,1.30200)--(9.62500,1.48650)--(9.73000,1.67100)--(9.83500,1.85550)--(9.94000,2.04000)--(10.04500,2.22450)--(10.15000,2.40900)--(10.25500,2.59350);
\draw[entrycolor,line width=1.4pt] (10.25500,2.59350)--(10.36000,2.77800);
\draw[slow] (10.36000,2.77800)--(10.46500,2.85959)--(10.57000,2.83103)--(10.67500,2.83085)--(10.78000,2.93351)--(10.88500,2.98129)--(10.99000,2.94037)--(11.09500,2.97719)--(11.20000,3.08267)--(11.30500,3.09272)--(11.41000,3.05870)--(11.51500,3.13147)--(11.62000,3.22079)--(11.72500,3.19906)--(11.83000,3.18939)--(11.93500,3.28767)--(12.04000,3.34606)--(12.14500,3.30679)--(12.25000,3.33256)--(12.35500,3.43925)--(12.46000,3.46003)--(12.56500,3.42204)--(12.67000,3.48510)--(12.77500,3.58089)--(12.88000,3.56720)--(12.98500,3.54899)--(13.09000,3.64143)--(13.19500,3.70987)--(13.30000,3.67385)--(13.40500,3.68874)--(13.51000,3.79498)--(13.61500,3.82676)--(13.72000,3.78639)--(13.82500,3.83906)--(13.93000,3.93991)--(14.03500,3.93530)--(14.14000,3.90969)--(14.24500,3.99495)--(14.35000,4.07265)--(14.45500,4.04142)--(14.56000,4.04582)--(14.66500,4.14997)--(14.77000,4.19278)--(14.87500,4.15165)--(14.98000,4.19349)--(15.08500,4.29789);
\fill[ink!40] (9.10000,0.56400) circle (.85pt);
\fill[ink!40] (9.20500,0.74850) circle (.85pt);
\fill[ink!40] (9.31000,0.93300) circle (.85pt);
\fill[ink!40] (9.41500,1.11750) circle (.85pt);
\fill[ink!40] (9.52000,1.30200) circle (.85pt);
\fill[ink!40] (9.62500,1.48650) circle (.85pt);
\fill[ink!40] (9.73000,1.67100) circle (.85pt);
\fill[ink!40] (9.83500,1.85550) circle (.85pt);
\fill[ink!40] (9.94000,2.04000) circle (.85pt);
\fill[ink!40] (10.04500,2.22450) circle (.85pt);
\fill[ink!40] (10.15000,2.40900) circle (.85pt);
\fill[entrycolor] (10.25500,2.59350) circle (.85pt);
\fill[slowcolor] (10.36000,2.77800) circle (.85pt);
\fill[slowcolor] (10.46500,2.85959) circle (.85pt);
\fill[slowcolor] (10.57000,2.83103) circle (.85pt);
\fill[slowcolor] (10.67500,2.83085) circle (.85pt);
\fill[slowcolor] (10.78000,2.93351) circle (.85pt);
\fill[slowcolor] (10.88500,2.98129) circle (.85pt);
\fill[slowcolor] (10.99000,2.94037) circle (.85pt);
\fill[slowcolor] (11.09500,2.97719) circle (.85pt);
\fill[slowcolor] (11.20000,3.08267) circle (.85pt);
\fill[slowcolor] (11.30500,3.09272) circle (.85pt);
\fill[slowcolor] (11.41000,3.05870) circle (.85pt);
\fill[slowcolor] (11.51500,3.13147) circle (.85pt);
\fill[slowcolor] (11.62000,3.22079) circle (.85pt);
\fill[slowcolor] (11.72500,3.19906) circle (.85pt);
\fill[slowcolor] (11.83000,3.18939) circle (.85pt);
\fill[slowcolor] (11.93500,3.28767) circle (.85pt);
\fill[slowcolor] (12.04000,3.34606) circle (.85pt);
\fill[slowcolor] (12.14500,3.30679) circle (.85pt);
\fill[slowcolor] (12.25000,3.33256) circle (.85pt);
\fill[slowcolor] (12.35500,3.43925) circle (.85pt);
\fill[slowcolor] (12.46000,3.46003) circle (.85pt);
\fill[slowcolor] (12.56500,3.42204) circle (.85pt);
\fill[slowcolor] (12.67000,3.48510) circle (.85pt);
\fill[slowcolor] (12.77500,3.58089) circle (.85pt);
\fill[slowcolor] (12.88000,3.56720) circle (.85pt);
\fill[slowcolor] (12.98500,3.54899) circle (.85pt);
\fill[slowcolor] (13.09000,3.64143) circle (.85pt);
\fill[slowcolor] (13.19500,3.70987) circle (.85pt);
\fill[slowcolor] (13.30000,3.67385) circle (.85pt);
\fill[slowcolor] (13.40500,3.68874) circle (.85pt);
\fill[slowcolor] (13.51000,3.79498) circle (.85pt);
\fill[slowcolor] (13.61500,3.82676) circle (.85pt);
\fill[slowcolor] (13.72000,3.78639) circle (.85pt);
\fill[slowcolor] (13.82500,3.83906) circle (.85pt);
\fill[slowcolor] (13.93000,3.93991) circle (.85pt);
\fill[slowcolor] (14.03500,3.93530) circle (.85pt);
\fill[slowcolor] (14.14000,3.90969) circle (.85pt);
\fill[slowcolor] (14.24500,3.99495) circle (.85pt);
\fill[slowcolor] (14.35000,4.07265) circle (.85pt);
\fill[slowcolor] (14.45500,4.04142) circle (.85pt);
\fill[slowcolor] (14.56000,4.04582) circle (.85pt);
\fill[slowcolor] (14.66500,4.14997) circle (.85pt);
\fill[slowcolor] (14.77000,4.19278) circle (.85pt);
\fill[slowcolor] (14.87500,4.15165) circle (.85pt);
\fill[slowcolor] (14.98000,4.19349) circle (.85pt);
\fill[slowcolor] (15.08500,4.29789) circle (.85pt);
\draw[fill=white,draw=entrycolor,line width=.9pt] (10.25500,2.59350) circle (2pt);
\draw[fill=white,draw=slowcolor,line width=.9pt] (10.36000,2.77800) circle (2pt);
\draw[fill=white,draw=slowcolor,line width=.9pt] (15.08500,4.29789) circle (2pt);
\draw[ink!35,dashed,line width=.5pt] (10.36000,2.77800)--(10.36000,0.40000);
\node[below=3pt] at (10.36000,0.40000) {$s+1$};
\draw[ink!35,dashed,line width=.5pt] (15.08500,4.29789)--(15.08500,0.40000);
\node[below=3pt] at (15.08500,0.40000) {$\tau$};
\draw[leader] (10.25500,2.59350)--(8.95000,1.75000);
\node[anchor=east,text=entrycolor] at (8.95000,1.75000) {$n=s$};
\node[anchor=west,text=slowcolor] at (11.65000,2.85000) {$|p_{n+1}-p_n|\leq\eta_j$};
\draw[dim] (15.83000,2.77800)--(15.83000,4.29500);
\node[anchor=west,align=left] at (15.95000,3.53650) {$\geq R_j-3C$};
\draw[decorate,decoration={brace,mirror,amplitude=4pt},draw=slowcolor!75] (10.36000,-0.25000)--(15.08500,-0.25000);
\node[below=6pt,text=slowcolor] at (12.72250,-0.25000) {$\tau-s-1$ small steps};
\end{tikzpicture}
\endgroup
\caption{Slow passage through an annulus in Lemma~\ref{lem:slow-annuli}.}
\label{fig:slow-passage}
\end{figure}

\FloatBarrier

\begin{lemma}\label{lem:slow-annuli}
Let $(p_n)_{n\geq0}$ be a sequence in $\mathbb R^d$ with
$|p_{n+1}-p_n|\leq C$ for some $C>0$.
Suppose that $R_j\to\infty$ and $\eta_j\geq0$, $\eta_j\to0$, and that
\[
 R_j+C\leq|p_n|\leq2R_j-C
 \quad\Longrightarrow\quad
 |p_{n+1}-p_n|\leq\eta_j
\]
for every $j$ and $n$. Then
\[
 \liminf_{n\to\infty}\frac{|p_n|}{n}=0.
\]
\end{lemma}

\begin{proof}
The conclusion is immediate if $(p_n)$ is bounded. Otherwise, fix
$j$ large enough that $R_j>3C$ and $|p_0|<R_j+C$.
Let
\[
 \tau=\min\{n\geq0:|p_n|\geq2R_j-C\},\qquad
 s=\max\{0\leq n<\tau:|p_n|\leq R_j+C\}.
\]
Unboundedness gives $\tau<\infty$, and the set defining $s$ contains
$0$. The global step bound gives
$|p_{s+1}|\leq R_j+2C<2R_j-C$, so $s+1<\tau$.
For $s+1\leq n<\tau$, the points lie in the specified annulus.
Consequently,
\[
 R_j-3C\leq|p_\tau|-|p_{s+1}|
 \leq\sum_{n=s+1}^{\tau-1}|p_{n+1}-p_n|
 \leq(\tau-s-1)\eta_j\leq\tau\eta_j.
\]
The first-arrival property gives $|p_\tau|\leq2R_j$. Thus
\begin{equation}\label{eq:slowdown-ratio}
 \frac{|p_\tau|}{\tau}
 \leq\frac{2R_j}{\tau}
 \leq\frac{2R_j}{R_j-3C}\eta_j\longrightarrow0.
\end{equation}
These arrival times tend to infinity, since
$|p_\tau|\geq2R_j-C\to\infty$ and
$|p_n|\leq|p_0|+nC$.
\end{proof}

\Needspace{10cm}
\subsection{Invariant graphs in dimension one}\label{subsec:invariant-graphs}

We state the translated-curve theorem of R{\"u}ssmann~\cite{Russmann1970}
for perturbations of $(x,p)\mapsto(x+p,p)$; see also
\cite[Section~6.1, Theorem~6.1]{KatokKrikorian2022}.

\begin{theorem}\label{thm:translated-curves}
There exist an integer $k_0\geq2$ and $\delta_0>0$ with the following
property. Let
\[
 G:\mathbb T\times(e^{-1},1)\longrightarrow\mathbb T\times\mathbb R
\]
be a smooth diffeomorphism onto its image. Suppose that $G$ has a lift
$\widetilde G$ satisfying
$\widetilde G(x+1,p)=\widetilde G(x,p)+(1,0)$ and
\[
 \|\widetilde G-\widetilde G_0\|_{C^{k_0}(\mathbb T\times(e^{-1},1))}
 <\delta_0,
 \qquad \widetilde G_0(x,p)=(x+p,p).
\]
Then there exist
$\gamma\in C^1\bigl(\mathbb T,(e^{-3/4},e^{-1/4})\bigr)$ and
$c\in\mathbb R$ such that the graph
$\Gamma=\{(x,\gamma(x)):x\in\mathbb T\}$ satisfies
\[
 G(\Gamma)=\{(x,\gamma(x)+c):x\in\mathbb T\}.
\]
Such a graph is called a \emph{translated graph}.
\end{theorem}

\begin{lemma}\label{lem:bnd-circle}
There exist an integer $k_*\geq2$ and $\delta_*>0$ such that, for every
$B\in C^\infty(\mathbb T)$ satisfying
$\lVert B\rVert_{C^{k_*}}<\delta_*$, the map
\begin{equation}\label{eq:bnd-residual-map}
 G_B(x,p)=\bigl(x+p+B'(x),\ p+B'(x)\bigr)
\end{equation}
has an invariant $C^1$ graph $p=\gamma(x)$ with $0<\gamma(x)<1$.
\end{lemma}

\begin{proof}
Let $k_0$ and $\delta_0$ be given by
Theorem~\ref{thm:translated-curves}. Set
$k_*=k_0+1$ and $\delta_*=\delta_0/\sqrt2$.
The map $G_B$ is a smooth diffeomorphism. Its lift satisfies
\[
 \widetilde G_B-\widetilde G_0=(B',B'),\qquad
 \|\widetilde G_B-\widetilde G_0\|_{C^{k_0}}
 \leq\sqrt2\|B\|_{C^{k_*}}<\delta_0.
\]
Theorem~\ref{thm:translated-curves} gives a translated $C^1$ graph
$\Gamma=\{(x,\gamma(x)):x\in\mathbb T\}$ contained in
$\mathbb T\times(e^{-3/4},e^{-1/4})$. Let $c$ be its translation.

Write $(X,P)=G_B(x,p)$. Direct calculation gives
\[
 P\,dX-p\,dx=d\bigl(B(x)+\tfrac12P^2\bigr).
\]
The induced circle map has degree one. Hence
\[
 0=\int_{G_B\Gamma}p\,dx-\int_\Gamma p\,dx
  =\int_0^1\bigl(\gamma(x)+c-\gamma(x)\bigr)\,dx=c.
\]
Thus $\Gamma$ is invariant.
\end{proof}

\begin{lemma}\label{lem:local-invariant-graph}
Let $S\in C^2(\mathbb R^2)$ satisfy diagonal periodicity and the
uniform twist condition, and let $U,V\in C^\infty(\mathbb T)$.
Suppose that, for some integer $a$,
\[
 S(x,y)=\frac12(y-x)^2+U(x)
 \qquad\text{whenever }a<y-x<a+1.
\]
Define $G_{U-V}:\mathbb T\times\mathbb R\to\mathbb T\times\mathbb R$ by
\[
 G_{U-V}(x,p)=\bigl(x+p+U'(x)-V'(x),\ p+U'(x)-V'(x)\bigr).
\]
If $G_{U-V}$ has an invariant graph $p=\gamma(x)$ with
$0<\gamma(x)<1$, then the map $F_V$ generated by $S-V(x)$
has the invariant graph
\[
 \Gamma_a=\{(x,a+\gamma(x)):x\in\mathbb T\}.
\]
\end{lemma}

\begin{proof}
Put $B=U-V$. For a lift $x$ and $p=a+\gamma(x)$, set
\[
 P=p+B'(x),\qquad Y=x+P.
\]
Then
\[
 \gamma(x)+B'(x)=P-a,\qquad
 x+\gamma(x)+B'(x)=Y-a.
\]
By the definition of $G_B$,
\[
 \begin{aligned}
 G_B(x\bmod1,\gamma(x))
 &=\bigl((x+\gamma(x)+B'(x))\bmod1,\ \gamma(x)+B'(x)\bigr)\\
 &=\bigl((Y-a)\bmod1,\ P-a\bigr).
 \end{aligned}
\]
Since $a\in\mathbb Z$, we have $(Y-a)\bmod1=Y\bmod1$.
The graph $p=\gamma(x)$ is invariant under $G_B$, so
\[
 P-a=\gamma(Y\bmod1),\qquad
 P=a+\gamma(Y\bmod1)\in(a,a+1).
\]
The last inclusion follows from $0<\gamma<1$.
Thus the displacement $Y-x=P$ lies where the stated formula for
$S$ holds. The generating equations of $S-V$ there are
$p=(Y-x)-B'(x)$ and $P=Y-x$, so the twist condition identifies
this step with $F_V$. Thus, on $\mathbb T\times\mathbb R$,
\[
 F_V(x,a+\gamma(x))=G_B(x,\gamma(x))+(0,a).
\]
Since $G_B$ maps its graph onto itself, $F_V(\Gamma_a)=\Gamma_a$.
\end{proof}

\begin{lemma}\label{lem:graph-barriers}
Let $F$ be a homeomorphism of $\mathbb T\times\mathbb R$.
Suppose that it has invariant continuous
graphs $p=\gamma_k^\pm(x)$ such that
\[
 \min_{x\in\mathbb T}\gamma_k^+(x)\longrightarrow+\infty,
 \qquad
 \max_{x\in\mathbb T}\gamma_k^-(x)\longrightarrow-\infty.
\]
Then every orbit of $F$ has bounded momentum.
\end{lemma}

\begin{proof}
Fix $z=(x_0,p_0)$. Choose invariant graphs lying entirely above and
below the level $p=p_0$. The region between them is the only connected
component of their complement with compact closure. Since $F$ is a
homeomorphism and preserves both graphs, it maps this component and
its closure onto themselves. Thus every iterate $F^n(z)$,
$n\in\mathbb Z$, lies in this compact annulus.
\end{proof}

\section{Small potentials with given momentum drift}
\label{sec:acceleration}

We apply the common drift ball from Lemma~\ref{lem:exclude-and-cancel}.

\begin{proof}[Proof of Theorem~\ref{thm:acceleration}]
Fix $\varepsilon>0$ and $0<\alpha<1$.
Apply Lemma~\ref{lem:potential-correction} with the constant $K$
from Assumption~\ref{ass:generating-function}.
Let $W,u,I_0,I_1,r$ be the resulting objects, and set $V=-W$.
Then $\|V\|_{C^{1,\alpha}}<\varepsilon$.

The functions $Q_{p^*}$ from Lemma~\ref{lem:potential-correction}
satisfy the hypotheses of Lemma~\ref{lem:exclude-and-cancel}.
The latter gives the required
orbit for each $|p^*|\leq r$, with the uniform bound
in~\eqref{eq:local-drift}.
The potential $V$ and the radius $r$ were chosen before $p^*$.
This proves~\eqref{eq:bounded-drift-error} for the whole closed ball.
\end{proof}

\begin{remark}
Although $\|V\|_{C^{1,\alpha}}$ is small,
\eqref{eq:negative-strip-potential} gives
$\|D^2V\|_\infty\geq2K+1$.
\end{remark}

\section{Persistence and generic linear momentum escape}\label{sec:persistence}

To prove Theorem~\ref{thm:robust-drift}, we absorb the given potential
into the generating function and apply Lemma~\ref{lem:drift-persistence}.

\begin{proof}[Proof of Theorem~\ref{thm:robust-drift}]
Fix $U\in C^\infty(\mathbb T^d)$, $\varepsilon>0$, and
$0<\alpha<1$. Set
\[
 S_U(x,y)=S(x,y)-U(x).
\]
The function $S_U$ satisfies Assumption~\ref{ass:generating-function}
with twist constant
$\lambda$ and bound $K+\|D^2U\|_\infty$.

Apply Lemma~\ref{lem:potential-correction} with $K$ replaced by
$K+\|D^2U\|_\infty$.
Let $W,u,I_0,I_1,r$ be the resulting objects, and set $V=U-W$.
Then
\[
 \|V-U\|_{C^{1,\alpha}}=\|W\|_{C^{1,\alpha}}<\varepsilon.
\]
The functions $Q_{p^*}$ from Lemma~\ref{lem:potential-correction}
satisfy the hypotheses of Lemma~\ref{lem:drift-persistence} for $S_U$ and $W$.
Apply that lemma to obtain $\delta>0$ and the bound
\[
 C=\frac{2\sqrt d\,\bigl(K+\|D^2U\|_\infty\bigr)^2}{\lambda}.
\]

Fix $E\in C^\infty(\mathbb T^d)$ with $\|E\|_\infty<\delta$. Since
\[
 S_U(x,y)-(-W+E)(x)=S(x,y)-(V+E)(x),
\]
the map supplied by Lemma~\ref{lem:drift-persistence} is $F_{V+E}$.
For every $|p^*|\leq r$, the lemma therefore gives a point $z_{E,p^*}$
satisfying~\eqref{eq:persistent-drift}.
The radius $r$ and the bound $C$ are independent of $E$.
\end{proof}

\begin{proof}[Proof of Corollary~\ref{cor:open-dense}]
Fix $0<\alpha<1$. For each $U\in C^\infty(\mathbb T^d)$ and
$\varepsilon>0$, Theorem~\ref{thm:robust-drift} gives $V$ and
$\delta,r,C>0$. Since $\|E\|_\infty\leq\|E\|_{C^{1,\alpha}}$,
every smooth potential in the $C^{1,\alpha}$ ball of radius $\delta$
about $V$ realizes every drift $|p^*|\leq r$.
The union of these balls is open in the relative $C^{1,\alpha}$ topology.
It is dense because $\|V-U\|_{C^{1,\alpha}}<\varepsilon$ for every
$U$ and $\varepsilon>0$.
\end{proof}

\section{An obstruction to linear momentum growth}
\label{sec:slowdown}

We construct $S$ from the sequence $(U_k)$ in
Lemma~\ref{lem:repeated-potentials}. For each $V$ with
$\|V\|_{C^2}\leq\varepsilon_0$, we choose a subsequence that converges
to $V$ in $C^2$. The differences $DU_k-DV$ control the momentum
increments on the corresponding annuli.
Lemma~\ref{lem:slow-annuli} then gives~\eqref{eq:zero-liminf}.

\begin{proof}[Proof of Theorem~\ref{thm:slowdown}]
Choose $\chi\in C_c^\infty(\mathbb R^d,[0,1])$ such that
\[
 \operatorname{supp}\chi\subset\{v:1/2<|v|<3\},
 \qquad \chi(v)=1\quad\text{if }1\leq|v|\leq2.
\]
Set $\chi_k(v)=\chi(6^{-k}v)$ for $k\geq1$.
Choose $\varepsilon_0>0$ so small that
\[
 2\varepsilon_0\|\chi\|_{C^2}\left(\frac16+\frac1{36}\right)
 \leq\frac12.
\]
Let $(U_k)$ be the sequence in Lemma~\ref{lem:repeated-potentials}
for this $\varepsilon_0$. Define
\begin{equation}\label{eq:radial-catalogue-S}
 S(x,y)=\frac12|y-x|^2+
       \sum_{k\geq1}\chi_k(y-x)U_k(x).
\end{equation}
The supports of the $\chi_k$ are pairwise disjoint and locally finite.
Moreover,
\[
 \|D\chi_k\|_\infty\leq6^{-k}\|\chi\|_{C^2},\qquad
 \|D^2\chi_k\|_\infty\leq6^{-2k}\|\chi\|_{C^2}.
\]
Lemma~\ref{lem:localized-generator} therefore shows that $S$ is smooth
and diagonally periodic, has bounded Hessian, and satisfies
\eqref{eq:shell-twist-bound}. In particular, it satisfies
Assumption~\ref{ass:generating-function}.

Fix $V\in C^\infty(\mathbb T^d)$ with
$\|V\|_{C^2}\leq\varepsilon_0$. By
\eqref{eq:catalogue-generating-equations}, every step of its map
satisfies
\begin{equation}\label{eq:catalogue-step-bound}
 |P-p|\leq3\varepsilon_0.
\end{equation}
This bound is independent of $V$ and of the phase point.

Lemma~\ref{lem:repeated-potentials} gives a strictly increasing
sequence $(k_j)$ such that
\begin{equation}\label{eq:matching-shells}
 \|U_{k_j}-V\|_{C^2}<\frac1j.
\end{equation}
If
\begin{equation}\label{eq:input-momentum-annulus}
 6^{k_j}+3\varepsilon_0\leq|p|\leq2\cdot6^{k_j}-3\varepsilon_0,
\end{equation}
then \eqref{eq:catalogue-step-bound} gives
$6^{k_j}\leq|P|\leq2\cdot6^{k_j}$. Thus $\chi_{k_j}(P)=1$ and
$D\chi_{k_j}(P)=0$.
Differentiating~\eqref{eq:radial-catalogue-S} in the second variable gives
\[
 D_2S(x,y)=y-x+\sum_{k\geq1}U_k(x)D\chi_k(y-x).
\]
Since the supports are disjoint, $D\chi_k(P)=0$ for $k\ne k_j$.
Substituting $x=y-P$ gives
\[
 D_2S(y-P,y)=P+U_{k_j}(y-P)D\chi_{k_j}(P)=P.
\]
By Lemma~\ref{lem:global-map}, the equation $D_2S(x,y)=P$ determines
$x$ uniquely for fixed $(y,P)$. Hence $x=y-P$ and
$\chi_{k_j}(y-x)=1$.
The generating equations~\eqref{eq:catalogue-generating-equations}
and the bound~\eqref{eq:matching-shells} now give
\begin{equation}\label{eq:matching-step}
 P-p=DU_{k_j}(x)-DV(x),\qquad |P-p|<\frac1j.
\end{equation}
For every initial point $z$, the sequence
$p_n=\pi_pF_V^n(z)$ satisfies the hypotheses of
Lemma~\ref{lem:slow-annuli}, with
\[
 R_j=6^{k_j},\qquad \eta_j=\frac1j,\qquad C=3\varepsilon_0.
\]
Therefore
\[
 \liminf_{n\to\infty}\frac{|\pi_pF_V^n(z)|}{n}=0.\qedhere
\]
\end{proof}

\section{Bounded momentum in dimension one}\label{sec:bounded}

We now place the repeated potentials at large positive and negative
integer displacements. Lemmas~\ref{lem:bnd-circle}
and~\ref{lem:local-invariant-graph} give invariant graphs in these
regions. Lemma~\ref{lem:graph-barriers} then bounds every orbit.

\begin{proof}[Proof of Theorem~\ref{thm:bounded}]
Choose $\beta\in C_c^\infty((-3,3))$ such that
$0\leq\beta\leq1$ and $\beta=1$ on $[-2,2]$.
Choose $\varepsilon_0>0$ so small that
\begin{equation}\label{eq:bnd-budget}
 4\varepsilon_0\|\beta\|_{C^2}\leq\tfrac12.
\end{equation}
Let $(U_k)$ be the sequence in Lemma~\ref{lem:repeated-potentials}
for $d=1$ and this $\varepsilon_0$.
For $k\geq1$, set
\[
 \chi_k(v)=\beta(v-6k)+\beta(v+6k).
\]
The supports are pairwise disjoint and locally finite, and
$\|\chi_k\|_{C^2}\leq\|\beta\|_{C^2}$. Define
\begin{equation}\label{eq:bnd-background}
 S(x,y)=\frac12(y-x)^2+\sum_{k\geq1}\chi_k(y-x)U_k(x).
\end{equation}
Lemma~\ref{lem:localized-generator} gives smoothness, diagonal
periodicity, and a bounded Hessian. It also gives
\[
 |D_{12}S+1|
 \leq4\varepsilon_0\|\beta\|_{C^2}\leq\tfrac12.
\]
In particular, $S$ satisfies Assumption~\ref{ass:generating-function}.

Fix $V\in C^\infty(\mathbb T)$ with
$\|V\|_{C^2}\leq\varepsilon_0$.
Let $k_*$ and $\delta_*$ be given by Lemma~\ref{lem:bnd-circle}.
Lemma~\ref{lem:repeated-potentials} gives a term $U_j$ such that
\[
 \|U_j-V\|_{C^{k_*}}<\delta_*.
\]
Let $p=\gamma(x)$ be the invariant graph of $G_{U_j-V}$ given by
Lemma~\ref{lem:bnd-circle}. Whenever $U_k=U_j$ and
$a\in\{-6k,6k\}$, we have $\chi_k(y-x)=1$ for
$a<y-x<a+1$. On this region,
\[
 S(x,y)=\frac12(y-x)^2+U_j(x).
\]
Lemma~\ref{lem:local-invariant-graph} therefore gives the invariant
graphs
\begin{equation}\label{eq:bnd-translated-circles}
 \Gamma_a=\{(x,a+\gamma(x)):x\in\mathbb T\},
 \qquad a\in\{-6k,6k\}.
\end{equation}
Since $U_j$ occurs infinitely often, these graphs have arbitrarily
large positive and negative momentum.
Lemma~\ref{lem:graph-barriers} now implies
\[
 \sup_{n\in\mathbb Z}|\pi_pF_V^n(z)|<\infty
\]
for every $z\in\mathbb T\times\mathbb R$.
\end{proof}

\begin{remark}
The sequence has a uniform $C^2$ bound, while its $C^{k_*}$
approximation property is used separately for each smooth potential.
Thus the proof does not require a $C^2$ version of the invariant-curve
theorem.
\end{remark}

\Needspace{10cm}
\section*{Acknowledgments}
The author thanks Yujia An for helpful discussions.

\section*{AI-use declaration}
The author used ChatGPT and Codex (OpenAI) to assist with language
editing, discussions of proofs, and LaTeX preparation. The author
is responsible for the content of the manuscript.

\section*{Data Availability}
The manuscript has no associated data.

\section*{Declarations}
\noindent\textbf{Conflicts of interest.}
The author declares no conflict of interest.

\bibliographystyle{amsplain}
\bibliography{reference}
\end{document}